\documentclass{amsart}
\usepackage{amssymb,latexsym,enumerate,amscd}

\swapnumbers
\theoremstyle{plain}
\newtheorem{theorem}{Theorem}

\newtheorem{corollary}[theorem]{Corollary}

\theoremstyle{definition}
\newtheorem{definition}[theorem]{Definition}

\newtheorem{Thm}[theorem]{Theorem}
\newtheorem{Lem}[theorem]{Lemma}
\theoremstyle{definition}
\newtheorem*{Not*}{Notation and Terminology}
\DeclareMathOperator{\Tor}{Tor}
\DeclareMathOperator{\Ext}{Ext}

\newcommand{\complexs}{\mathbb{C}}

\newcommand{\integers}{\mathbb{Z}}

\DeclareMathOperator{\id}{id}

\newcommand{\tensor}{\otimes}

\newcommand{\innerprod}[1]{\langle #1 \rangle}

{\catcode`@=11\global\let\c@equation=\c@theorem}

\allowdisplaybreaks[2]

\begin{document}

\title[The Ore condition and lamplighter group]{The
Ore condition, affiliated operators,
and the lamplighter group}

\author[P. A. Linnell]{Peter A. Linnell}
\address{Department of Mathematics, Virginia Tech\\
Blacksburg, VA 24061-0123, USA}
\email{linnell@math.vt.edu}
\urladdr{http://www.math.vt.edu/people/linnell/}
\thanks{The first author was supported in part by the
Sonderforschungsbereich in M{\"u}nster}

\author[W. L{\"u}ck]{Wolfgang L{\"u}ck}
\address{FB Mathematik, Universit{\"a}t M{\"u}nster\\
Einsteinstr. 62, D-48149 M{\"u}nster, Germany}
\email{lueck@math.uni-muenster.de}
\urladdr{http://wwwmath.uni-muenster.de/math/u/lueck/}

\author[T. Schick]{Thomas Schick}
\address{Mathematisches Institiut, Universit{\"a}t G{\"o}ttingen\\
Bunsensrt. 3-5, D-37073 G{\"o}ttingen, Germany}
\email{schick@uni-math.gwdg.de}
\urladdr{http://www.uni-math.gwdg.de/schick/}

\begin{abstract}
  Let $G= \integers/2\integers \wr \integers$ be the so called
  lamplighter group and $k$ a commutative ring. We show that $kG$ does
  not have a classical ring of quotients (i.e.~does not satisfy the
  Ore condition). This answers a
Kourovka notebook problem.
Assume that $kG$ is contained in a ring $R$ in
  which the element $1-x$ is invertible, with $x$ a generator of
  $\integers\subset G$. Then $R$ is not flat over $kG$. If
  $k=\complexs$, this applies
  in particular to the algebra $\mathcal{U}G$ of unbounded operators
  affiliated to the group von Neumann algebra of $G$.

We present two proofs of these results. The second one is due to
Warren Dicks, who, having seen our argument, found a much simpler and
more elementary proof, which at the same time yielded a more general
result than we had originally proved. Nevertheless, we present both
proofs here, in the hope
that the original arguments might be of use in some other context not
yet known to us.

\end{abstract}

\keywords{Ore ring, affiliated operators, flat, lamplighter group,
Fox calculus}

\subjclass[2000]{Primary: 16U20; Secondary: 46L99}
%requires amsart.cls at least 2.0

\date{%CHANGE!!
August 4, 2003}

\maketitle

\section{Notation and Terminology}

Let $A$ be a group.  Then $A \wr \mathbb {Z}$ indicates the
Wreath product with base group
$B=\bigoplus_{i=-\infty}^{\infty} A_i$, where $A_i = A$ for
all $i$.  Thus $A \wr
\mathbb {Z}$ is isomorphic to the split extension $B \rtimes
\mathbb {Z}$ where if $x$ is a generator for $\mathbb {Z}$,
then $x^iA_0x^{-i} = A_{i}$.  Also we identify $A$ with
$A_0$.  In the case $A = \mathbb {Z}/2\mathbb {Z}$ above,
$A \wr \mathbb {Z}$ is often called the lamplighter group.

Let $kG$ denote the group algebra of the group $G$ over the
field $k$, and let $\alpha \in kG$.  Write $\alpha =
\sum_{g\in G} a_gg$ where $a_g \in k$.  Then the support of
$\alpha$ is $\{g \in G \mid a_g \ne 0\}$, a finite
subset of $G$.

The augmentation
ideal of a group algebra will denoted by the small German letter
corresponding to the capital Latin letter used to name the group.
Thus if $k$ is a field and $G$ is a group,
then $\mathfrak{g}$ is the ideal of $kG$ which has $k$-basis
$\{g-1 \mid g \in G \setminus 1\}$.
For the purposes of this paper, it will always be clear over which
field we are working when considering augmentation ideals.

\begin{definition} \label{DOre}
  A ring $R$ satisfies the \emph{Ore condition} if for any $s,t\in R$
  with $s$ a non-zerodivisor, there are $x,y\in R$ with $x$ a non-zerodivisor 
such that $sy=tx$. Formally, this means that $s^{-1}t =
  yx^{-1}$, and the condition makes sure that a classical ring of
  quotients, inverting all non-zerodivisors of $R$, can be constructed.
\end{definition}

Then $R$ is an Ore ring means that $R$
satisfies the Ore condition.
Equivalently this means that if $S$ is the set of
non-zerodivisors of $R$, then given $r \in R$ and $s \in S$, we
can find $r_1 \in R$ and $s_1 \in S$ such that $rs_1 = sr_1$.
In this situation we can form the Ore localization $RS^{-1}$, which
consists of elements $rs^{-1}$ where $r \in R$ and $s \in S$.

The above definition is really the right Ore condition, though
for group rings the right and left Ore conditions are equivalent.

In this note, we study, which group rings satisfy the Ore condition. It is well
known that this fails for a non-abelian free group.

On the other hand, abelian groups evidently satisfy the Ore
condition. In this note we show that the lamplighter groups (and
relatives) do not satisfy it. Note, however, that these groups are
$2$-step solvable, i.e.~close relatives of abelian groups.

Let $G$ be a group, let $\mathcal {N}(G)$ denote the group von
Neumann algebra of $G$, let
$\mathcal {U}(G)$ denote the algebra of
unbounded linear operators affiliated to $\mathcal {N}(G)$, and let
$\mathcal{D}(G)$ denote the division closure of $\mathbb {C}G$ in
$\mathcal {U}(G)$.  For more information on these notions,
see \cite[\S 8 and \S 9]{Linnell98} and
\cite[\S 8 and \S 10]{Lueck01}.  In particular we have the inclusion
of $\mathbb{C}$-algebras
\[
\mathbb{C}G \subseteq \mathcal{D}(G) \subseteq \mathcal {U}(G)
\]
and it is natural to ask whether $\mathcal{D}(G)$ and
$\mathcal{U}(G)$ are flat over $\mathbb {C}G$.

We use the following well-known and easily verified statement
without further comment in this paper.  If $k$ is a field and
$g\in G$ has infinite order, then $1-g$ is a non-zerodivisor in $kG$,
and in the case $k= \mathbb {C}$ we also have that $1-g$ is
invertible in $\mathcal{D}(G)$.

If $H$ is the nonabelian free group of rank 2, then we have an
exact sequence of $\mathbb{C}H$-modules
$0 \to \mathbb {C}H^2 \to \mathbb {C}H \to \mathbb {C} \to 0$.
It was shown in \cite[Theorem 1.3]{Linnell93} (see also
\cite[Theorem 10.2]{Linnell98} and \cite[Theorem 10.19 and Lemma
10.39]{Lueck01})
that $\mathcal{D}(H)$ is a division ring, so when we apply
$\otimes_{\mathbb{C}H} \mathcal{D}(H)$ to this sequence, it
becomes $\mathcal{D}(H)^2\to \mathcal{D}(H)\to 0\to 0$, since $(1-x)$
is invertible in $D(H)$ for every element $1\ne x\in H$, but $(1-x)$
acts as the zero operator on $\mathbb{C}$.  By counting dimension,
we get from this (adding the kernel) a short exact sequence
$0 \to \mathcal{D}(H) \to \mathcal{D}(H)^2 \to
\mathcal{D}(H) \to 0$.  Suppose $Q$ is a ring containing
$\mathcal{D}(H)$ which is flat over $\mathbb{C}H$.  Then
applying $\otimes_{\mathcal{D}(H)} Q$
to the previous sequence, we obtain the exact sequence
$0 \to Q \to Q^2 \to Q \to 0$ which contradicts the
hypothesis that $Q$ is flat over $\mathbb{C}H$ (in the latter case we
would have obtained $0\to Q^2\to Q\to 0\to 0$).  In
particular if $G$ is a group containing $H$, then neither
$\mathcal{D}(G)$ nor $\mathcal{U}(G)$ is flat over
$\mathbb {C}H$.  Since
$\mathbb{C}G$ is a free $\mathbb {C}H$-module, we conclude that
neither $\mathcal{D}(G)$ nor $\mathcal {U}(G)$ is flat over
$\mathbb {C}G$.

To sum up the previous paragraph, neither $\mathcal {D}(G)$ nor
$\mathcal {U}(G)$ is flat over $\mathbb {C}G$ when $G$ contains
a nonabelian free group.
On the other hand it was proven in \cite[Theorem
9.1]{Reich98} (see also \cite[Theorem 10.84]{Lueck01}) that if $G$ is
an elementary amenable group which has a bound on the orders of its
finite subgroups, then $\mathcal{D}(G)$ and $\mathcal{U}(G)$ are flat
over $\mathbb{C}G$.  Furthermore it follows from \cite[Theorems
6.37 and 8.29]{Lueck01} that if $G$ is amenable (in particular if
$G$ is the lamplighter group), then at least $\mathcal{U}(G)$ is
``dimension flat" over $\mathbb {C}G$.

\section{Original results and proof}
\label{sec:orig-results-proof}

We shall prove

\begin{Thm} \label{Tnotflat}
Let $H\ne 1$ be a finite group and let $G$ be a group containing
$H \wr \mathbb {Z}$.  Then
neither $\mathcal {D}(G)$ nor $\mathcal {U}(G)$ is flat over $\mathbb
{C}G$.
\end{Thm}

Closely related to this question is the problem of when
the group algebra $kG$ of the group $G$ over the field $k$ is
an Ore ring (in other words does $kG$ have a classical ring of
quotients; see Definition \ref{DOre}).  Our next result
answers a Kourovka Notebook problem \cite[12.47]{MazurovKhukhro99},
which was proposed by the first author.  The problem there asks if
$kG$ has a classical quotient ring in the case $G= \mathbb {Z}_p \wr
\mathbb {Z}$ where $p$ is a prime.
\begin{Thm} \label{Tnotore}
Let $H\ne 1$ be a finite group, let $k$ be a field
and let $G$ be a group containing $H \wr \mathbb {Z}$.  Then
$kG$ is not an Ore ring.
\end{Thm}

\begin{Lem} \label{Lfingenproj}
Let $R$ be a subring of the ring $S$ and let $P$ be a projective
$R$-module.  If $P \otimes_R S$ is finitely generated as an
$S$-module, then $P$ is finitely generated.
\end{Lem}
\begin{proof}
Since $P$ is projective, there are $R$-modules $Q,F$ with $F$ free
such that $P \oplus Q = F$.  Let $\mathcal {E}$ be a basis for $F$.
Now $P \otimes_R S \oplus Q \otimes_R S = F \otimes_R S$ and since
$P \otimes_R S$ is finitely generated,
there exist $e_1, \dots, e_n \in \mathcal
{E}$ such that $P \otimes_R S \subseteq e_1 S + \dots + e_n S$.
We now see that every element $p$ of $P$ is
\begin{enumerate} [\normalfont (i)]
\item An $R$-linear combination of elements in $\mathcal {E}$
  ($\implies p\otimes 1= \sum_{e\in E} e \otimes r_e$).
\item An $S$-linear combination of $e_1, \dots, e_n$ ($p\otimes 1=
  \sum e_i\otimes s_i$).
\end{enumerate}
Set $E = e_1R + \dots + e_nR$.  Comparing coefficients, the above shows that
$P \subseteq E$ and it follows that we have the equation $P \oplus
(Q \cap E) = E$.
Therefore $P$ is a finitely generated $R$-module as required.
\end{proof}

\begin{Lem} \label{Lnotflat}
Let $H$ be a nontrivial finite group, let $k$ be a field with
characteristic which does not divide $|H|$, and let
$G = H \wr \mathbb {Z}$.  If $Q$
is a ring containing $kG$ such that $k \otimes_{kG} Q = 0$, then
$\Tor^{kG}(k,Q) \ne 0$.
\end{Lem}

\begin{proof}
Let $d$ denote the minimum
number of elements required to generate $G$.
Then we have exact sequences
\begin{gather*}
0 \longrightarrow \mathfrak{g} \longrightarrow kG \longrightarrow k
\longrightarrow 0 \\
0 \longrightarrow P \longrightarrow kG^d \longrightarrow \mathfrak
{g} \longrightarrow 0.
\end{gather*}
Suppose to the contrary $\Tor^{kG}(k,Q) = 0$.  Then the following
sequence is exact:
\begin{equation*}
  0 \longrightarrow \mathfrak{g} \otimes_{kG} Q \longrightarrow kG
\otimes_{kG} Q \longrightarrow k\otimes_{kG} Q = 0.
\end{equation*}
Since $kG$ has
homological dimension one \cite[p.~70 and Proposition 4.12]{Bieri81},
we also have
$0=\Tor_2^{kG}(k,Q) = \Tor_1^{kG}(\mathfrak{g},Q)$.  Hence we have another
exact sequence
\begin{equation}\label{eq:s1}
0 \longrightarrow P\otimes_{kG} Q  \longrightarrow kG^d
\otimes_{kG} Q
\longrightarrow \mathfrak {g}\otimes_{kG} Q  \longrightarrow 0.
\end{equation}
We rewrite this to get the exact sequence
\[
0 \longrightarrow P \otimes_{kG} Q \longrightarrow Q^d
\longrightarrow Q \longrightarrow 0
\]
and we conclude that $P
\otimes_{kG} Q$ is a finitely generated $Q$-module, which is
projective since the sequence \eqref{eq:s1} splits.
Now $kG$ has cohomological dimension $\le 2$ \cite[p.~70 and
Theorem 4.6 and Proposition 4.12]{Bieri81},
hence $P$ is a projective $kG$-module. (To see this, let $P'\to P$ be
a map from a projective $kG$-module $P'$ onto $P$.  Because
$\Ext^2_{kG}(k,Q)=0$, this extends to a map $P'\to kG^d\to P$.  Since
the image of the first arrow is $P$, this gives a split of the
injection $P\to kG^d$, i.e.\ $P$ is projective.) Therefore
$P$ is finitely generated by Lemma \ref{Lfingenproj}.  But it is well
known that $P\cong R/[R,R]\otimes k$ as $kG$-modules, where $1\to R\to
F\to G\to 1$ is an exact
sequence of groups and $F$ is a free group with $d$ generators
(compare the proof of \cite[(5.3) Theorem]{Brown82}).  Consequently,
$G$ is almost finitely presented over $k$ as defined in
\cite{BieriStrebel78} if $P$ is a finitely
generated $kG$-module.  But this is a contradiction to
\cite[Theorem~A or Theorem~C]{BieriStrebel78}
(here we use $H \ne 1$), where the structure of almost finitely
presented groups such as $G$ is determined.
\end{proof}

\begin{corollary} \label{Cnotflat}
Let $H$ be a nontrivial finite group, let $k$ be a field with
characteristic which does not divide $|H|$, and let $G$ be a group
containing $H \wr \mathbb{Z}$.  Let $x$ be a generator for $\mathbb
{Z}$ in $G$.  If $Q$ is a ring containing $kG$ and $1-x$ is invertible
in $Q$, then $Q$ is not flat over $kG$.
\end{corollary}
\begin{proof}
Set $L = H \wr \mathbb {Z}$.  Since $1-x$ is invertible in $Q$,
we have $\mathbb{C} \otimes_{\mathbb{C}L} Q = 0$, so
from Lemma \ref{Lnotflat} we deduce that $\Tor^{kL}(k,Q) \ne 0$.
Furthermore $kG$ is flat over $kL$, consequently
$\Tor^{kG}(k \otimes_{kL} kG,Q) \ne 0$ by \cite[p.~2]{Bieri81},
in particular $Q$ is not flat over $kG$ as required.
\end{proof}

\begin{proof}[Proof of Theorem \ref{Tnotflat}]
Set $L = H \wr \mathbb {Z}$ and
let $x$ be a generator for $\mathbb {Z}$ in $L$.
Then $1-x$ is
a non-zerodivisor in $\mathcal{N}(G)$ and therefore is invertible
in $\mathcal{U}(G)$, and hence also invertible in $\mathcal{D}(G)$.
The result now follows from Corollary \ref{Cnotflat}.
\end{proof}

\begin{Lem} \label{Lfirstp}
Let $p$ be a prime, let $k$ be a field of characteristic $p$,
let $A$ be a group of order $p$, let
$G = A \wr \mathbb {Z}$ with base group $B$,
let $a$ be a generator for $A$,
and let $x \in G$ be a generator for $\mathbb{Z}$.  Then
there does not exist $\alpha,\sigma \in kG$ with $\sigma
\notin \mathfrak{b}kG$ such that $(1-a)\sigma = (1-x)
\alpha$.
\end{Lem}
\begin{proof}
Suppose there does exist $\alpha$ and $\sigma$ as above.
Observe that $\alpha \in \mathfrak{b} kG$, since $\mathfrak{b}kG$ is
the kernel of the map $kG\to k\mathbb{Z}$ induced from the obvious
projection $G=A\wr \mathbb{Z} \to \mathbb{Z}$ mapping $x$ to $x$, and
since $1-x$ is not a
zerodivisor in $k\mathbb{Z}$.

Thus we may write $\sigma = \tau + \sum_i s_ix^i$ where
$s_i \in k$, $\tau \in \mathfrak{b}kG$ and not all the
$s_i$ are zero, and $\alpha = \sum_i x^i \alpha_i$ where
$\alpha_i \in \mathfrak{b}$.  Then the equation
$(1-a)\sigma = (1-x)\alpha$ taken mod $\mathfrak{b}^2kG$
yields
\[
\sum_i (1-a)s_i x^i = \sum_i (1-x) x^i \alpha_i.
\]
Set $b_i = 1 - x^{-i} a x^i$.
By equating the coefficients of $x^i$, we obtain
$s_i b_i = \alpha_i - \alpha_{i-1}$ mod $\mathfrak{b}^2kG$
for all $i$.  Since $\alpha_i \ne 0$ for only finitely many
$i$, we deduce that $\sum s_ib_i = 0$ mod $\mathfrak{b}^2$.
Also $\mathfrak{b}/\mathfrak{b}^2 \cong B \otimes k$ as $k$-vector
spaces via the map induced by $b-1 \mapsto b \otimes 1$
and the elements $x^{-i}ax^i \otimes 1$
are linearly independent in
$B \otimes k$, consequently the $b_i$ are linearly independent
over $k$ mod $\mathfrak {b}^2$.  We now have a contradiction
and the result follows.
\end{proof}

\begin{Lem} \label{Lzerodivisor}
Let $k$ be a field,
let $H$ be a locally finite subgroup of the group $G$,
and let $\alpha_1, \dots, \alpha_n \in \mathfrak{h}kG$.
Then there exists $\beta \in kH
\setminus 0$ such that $\beta \alpha_i = 0$ for all $i$.
\end{Lem}

\begin{proof}
Let $T$ be a right transversal for $H$ in $G$, so $G$ is the
disjoint union of $\{Ht \mid t \in T\}$.  For each $i$, we
may write $\alpha_i = \sum_{t \in T} \beta_{it}t$ where
$\beta_{it} \in \mathfrak{h}$.   
Let $B$ be the subgroup generated by the supports of the
$\beta_{it}$.  Then $B$ is a finitely generated subgroup
of $H$ and thus $B$ is a finite $p$-group.
Also $\beta_{it} \in \mathfrak{b}$ for all $i,t$.  Set $\beta
= \sum_{b\in B} b$.  Then $\beta \mathfrak{b} = 0$ and the result
follows.
\end{proof}

\begin{Lem} \label{Lnotore}
Let $p$ be a prime, let $k$ be a field of characteristic $p$,
let $A$ be a group of order $p$, and let $G$ be a group containing $A
\wr \mathbb {Z}$.  Then $kG$ does not satisfy the Ore condition.
\end{Lem}

\begin{proof}
Let $H = A \wr \mathbb {Z}$, which we may regard as a subgroup of
$G$, with $x\in H\subset G$ a generator of $\mathbb{Z}$.  Let $B$ be
the base group of $H$ and let $T$ be a right
transversal for $H$ in $G$, so $G$ is the disjoint union of $\{Ht
\mid t\in T\}$.  Note that $1-x$ is a non-zerodivisor in $kG$.
Let $a$ be a generator for $A$.
Suppose $(1-a)\sigma = (1-x)\alpha$ where $\alpha, \sigma \in kG$
and $\sigma$ is a non-zerodivisor in $kG$.  Then we may write $\alpha
= \sum_{t \in T} \alpha_tt$ and $\sigma = \sum_{t \in T} \sigma_tt$
with $\alpha_t, \sigma_t \in kH$, and then we have
\[
(1-a) \sigma_t = (1-x) \alpha_t
\]
for all $t \in T$.  If $\sigma_t \in \mathfrak{b}kH$ for all $t \in
T$, then by Lemma \ref{Lzerodivisor} we see that there exists $\beta
\in kB \setminus 0$ such that $\beta \sigma_t = 0$ for all $t$.  This
yields $\beta \sigma = 0$ which contradicts the hypothesis that
$\sigma$ is a non-zerodivisor.  Therefore we may assume that there
exists $s \in T$ such that $\sigma_s \notin \mathfrak{b} kH$.  But
now the equation
\[
(1-a) \sigma_s = (1-x) \alpha_s
\]
contradicts Lemma \ref{Lfirstp}, and the result follows.
\end{proof}

\begin{proof}[Proof of Theorem \ref{Tnotore}]
If the characteristic of $k$ is $p$ and divides $|H|$, then the result
follows from Lemma 
\ref{Lnotore}.  On the other hand if $p$ does not divide $|H|$,
we suppose that $kG$ satisfies the Ore condition.  Then $kG$ has a
classical ring of quotients $Q$. It is a well known fact that such a
classical ring of quotients is always flat over its base (compare
\cite[p.~57]{Stenstrom(1975)}). In
particular, $Q$ is flat over $kG$.  Let $x$ be
a generator for $\mathbb {Z}$ in $G$.  Since $1-x$ is
a non-zerodivisor in $kG$, we see that $1-x$ is invertible in $Q$.
We now have a
contradiction by Lemma \ref{Lnotflat} and the result follows.
\end{proof}

\section{Warren Dicks' proof}
\label{sec:warren-dicks-proof}

In this section, we give our account of Warren Dicks' proof and
generalization of the results presented in Section
\ref{sec:orig-results-proof}. All credit has to go to him, all
mistakes are ours.
More precisely, we prove the following theorem
(for elements $x,y$ in a group,
  we use the commutator convention $[x,y]=xyx^{-1}y^{-1}$):

\begin{theorem}\label{theo:main}
  Let $2 \le d \in \mathbb {Z}$, let
$G=\innerprod{a,x\mid a^d=1,\, [a,x^lax^{-l}]=1;\;l=1,2,\ldots}$ be the wreath
  product $\integers/ d\integers \wr \integers$, and let $k$ be a
  nonzero commutative ring with unit.   If $u,v \in kG$ are such that
$u(a-1)=v(x-1)$, then $u$ is a left zerodivisor in $kG$.
In particular, $kG$ does not satisfy the Ore condition.
\end{theorem}
\begin{proof}
  For the last statement, note that
$(x-1)$ is a non-zerodivisor in
  $kG$ because $x$ has infinite order.

Recall that any presentation $H=\innerprod{S\mid R}$ of a
  group $H$ gives rise to an 
  exact sequence of left $kH$-modules
  \begin{equation}\label{eq:seq}
    \bigoplus_{r\in R} kH \xrightarrow{F} \bigoplus_{s\in S} kH
    \xrightarrow{\alpha}  kH \xrightarrow{\epsilon} k \to 0.
  \end{equation}
  Here, $\epsilon$ is the augmentation map defined by $\epsilon(h) =
1$ for all $h \in H$, $\alpha$ maps
  $u\overline{s}\in \bigoplus_{s\in S} kH$ (with $u\in kH$ and
  $\overline{s}$ the canonical basis element
    corresponding to the generator $s\in S$) to $u(s-1)\in kH$, and the
    map $F$ is given by the Fox calculus, i.e.~$u\overline{r}$ (where
    $u\in kH$ and 
    $\overline{r}$ is the canonical basis element corresponding to the
    relator $r\in R$) is mapped to 
    \begin{equation*}
      \sum_{s\in S} u\frac{\partial r}{\partial s}
      \overline{s}.
    \end{equation*}
    If $r=s_{i_1}^{\epsilon_1}\dots s_{i_n}^{\epsilon_n}$ with $s_i\in
    S$ and $\epsilon_i\in\{-1,1\}$, then the Fox derivative is defined by
    \begin{equation*}
      \frac{\partial r}{\partial s} := \sum_{k=1}^n 
      s_{i_1}^{\epsilon_1}\dots s_{i_{k-1}}^{\epsilon_{k-1}}
      \frac{\partial s_{i_k}^{\epsilon_k}}{\partial s}.
    \end{equation*}
    Here $\partial s/\partial s=1$, $\partial s^{-1}/\partial s=
    -s^{-1}$ and $\partial t^{\epsilon}/\partial s=0$ if $s\ne t\in S$
and $\epsilon = \pm 1$.

    The above sequence can be considered as the cellular chain complex
    (with coefficients $k$)
    of the universal covering of the standard presentation CW-complex
    given by $\innerprod{S\mid R}$. Since this space is $2$-connected,
    its first homology vanishes and its zeroth homology is isomorphic
    to $k$ (by the augmentation), which implies that the sequence
    is indeed exact. An outline of the proof can be found in
    \cite[II.5 Exercise 3]{Brown82} or in \cite[IV.2, Exercises]{Brown82}

    Now we specialize to the group $G$.
Let us write $r_0 = a^d$ and $r_l = [a,x^lax^{-l}]$ for $l \ge 1$.
    Suppose $u,v\in kG$ with $u(a-1)=v(x-1)$. Then
$\alpha(u\overline{a} - v\overline{x})=0$.
    Exactness implies that there exists a positive integer $N$ and
$z_l \in kG$ $(0 \le l \le N)$ such
    that $F(\sum_l z_l \overline{r_l}) = u\overline{a}-
    v\overline{x}$. We want to prove that $u$ is a zerodivisor.
Therefore we are concerned only with the $\overline{a}$
    component of $F(\sum_{0 \le l \le N} z_l \overline{r_l})$. This means we
    first must compute $\partial r/\partial a$ for all the relators in
    our presentation of $G$. This is easily done:
    \begin{align}
      \frac{\partial a^d}{\partial a} &= 1 + a + \cdots + a^{d-1}\\
      \frac{\partial [a, x^lax^{-l}]}{\partial a} &= \frac{\partial(ax^lax^{-l} a^{-1}
        x^l a^{-1}x^{-l})}{\partial a}\\
      & = 1 + ax^l - ax^lax^{-l} a^{-1} -
      ax^lax^{-l} a^{-1} x^l a^{-1},
    \end{align}
    the latter for $l>0$. Using the fact that $x^lax^{-l}$ commutes
    with $a$ for each $l$, we can simplify:
    \begin{equation*}
      \frac{\partial [a, x^lax^{-1}]}{\partial a} = 1 + ax^l - x^lax^{-l}  -
      x^l = x^l(x^{-l}ax^l-1) - (x^lax^{-l}-1). 
    \end{equation*}
Since $u$ is the coefficient of $\overline{a}$
    in $F(\sum_{0 \le l \le N} z_l \overline{r_l})$ we see that
\[
u=
z_0(1+a+ \dots + a^{d-1}) + \sum_{n=1}^N z_n x^n(x^{-n}ax^n-1) -
z_n(x^nax^{-n}-1).
\]
Now let $C = \langle x^nax^{-n} \mid 1 \le |n| \le N\rangle$, a
finite subgroup of the base group $\bigoplus_{i=-\infty}^{\infty}
\mathbb{Z}/d\mathbb{Z}$.  Then
$C\langle a \rangle = C \times \langle a \rangle$.  Set $\gamma =
(1-a) \sum_{c\in C} c$.  Then $\gamma \ne 0$ and $\beta \gamma = 0$.
We conclude that $u$ is a left zerodivisor
in $kG$ and the result follows.
\end{proof}

\begin{corollary} \label{CDicks}
Let $2 \le d \in \mathbb {Z}$ and let $G$ be a group containing
$\mathbb {Z}/d \mathbb {Z} \wr \mathbb {Z}$, and let $x \in G$ be a
generator for $\mathbb {Z}$.  Let $k$ be a nonzero
commutative ring with unit and let $Q$ be a ring containing $kG$ such
that $1-x$ becomes invertible in $Q$.  Then $Q$ is not flat over $kG$.
\end{corollary}

\begin{proof}
Since $kG$ is free, hence flat, as left
$k[\mathbb {Z}/d \mathbb {Z} \wr \mathbb {Z}]$-module,
by \cite[p.~2]{Bieri81} we may assume that $G =
\mathbb {Z}/d \mathbb {Z} \wr \mathbb {Z}$.  Now
  tensor the exact sequence \eqref{eq:seq} over $kG$ with $Q$. Then
the resulting sequence will also be exact and
  $\overline{a} - (a-1)(x-1)^{-1}\overline{x}$ will be in the kernel of
  $\id_Q\tensor \alpha$.  Therefore
  $\overline{a} - (a-1)(x-1)^{-1}\overline{x}$ will be in the image of
$\id_Q\tensor F$.  However the proof of Theorem \ref{theo:main}
  shows that if $u\overline{a} - v\overline{x}$ is in the image of
  $\id_Q\tensor F$, then $u$ is a zerodivisor; the only change is
that we want $z_l \in Q$ rather than $z_l \in kG$.
  Since $u=1$ in this situation, which is not a zerodivisor, the
tensored sequence is not exact.
\end{proof}
In the case $k$ is a subfield of $\mathbb{C}$ and $Q= \mathcal{D}(G)$
or $\mathcal{U}(G)$, Corollary \ref{CDicks} tells us that if $G$
contains $\mathbb {Z}/d \mathbb {Z} \wr \mathbb {Z}$, then
$\mathcal{D}(G)$ and $\mathcal{U}(G)$ are not flat over $kG$.

\bibliographystyle{plain}
\bibliography{lamplighter}

\end{document}